\documentstyle{amsppt} 
\magnification=\magstep1 
\NoRunningHeads
\loadbold
\def\F{\Cal F}
\def\G{\Cal G}
\def\NN{\Bbb N}
\def\ep{\epsilon >0}
\def\e{\epsilon}
\def\l1{{\Cal L}_{E}^{1}}
\def\rl1{{\Cal L}_{\Bbb R}^{1}}
\def\dl1{{\Cal L}_{E^*}^\infty [E]}
\def\li{{\Cal L}_{E^*}^{\infty}}
\def\p1{{\Cal P}^{1}_{E}}
\def\pn{\text{Pettis }}
\def\qd{\hfill\ \qed}
\topmatter
\title
From weak to strong types of  
${\fam \cmbsyfam L}^{\boldkey 1}_{\boldkey E}$-convergence 
by the Bocce-criterion
\endtitle
\author
Erik J. Balder  \quad 
Maria Girardi$^*$   \quad
Vincent Jalby$^{**}$   
\endauthor
\address
Mathematical Institute, University of Utrecht, P.O. Box 80.010,
3508 TA Utrecht, Netherlands
\endaddress
\email balder\@math.ruu.nl \endemail
\address
Department of Mathematics, University of South Carolina, 
Columbia, South Carolina 29208, U.S.A.
\endaddress
\email girardi\@math.scarolina.edu \endemail
\address 
D\'epartement de Math\'ematiques, 
Universit\'e Montpellier II, 
34095 Montpellier Cedex 05, France
\endaddress   
\curraddr    
D\'epartement de Math\'ematiques, Universit\'e de Limoges,
87060 Limoges Cedex, France
\endcurraddr
\email jalby\@math.univ-montp2.fr \endemail
\thanks
$^*$ Supported in part by NSF grant  DMS-9204301 and DMS-9306460. 
\endthanks
\thanks
$^{**}$  Part of this 
paper is contained in the doctoral dissertation
of this author, who is 
grateful for the guidance of his advisor Professor C. Castaing.   
\endthanks
\date  
14 February 1994 
\enddate
\abstract   
Necessary and sufficient  oscillation conditions are given 
for a weakly convergent sequence  
(resp. relatively weakly compact set)  
in  the Bochner-Lebesgue  
space $\l1$ to be norm convergent 
(resp. relatively norm compact), 
thus extending the  known results  for $\rl1$.    
Similarly, 
necessary and sufficient oscillation conditions are given  
to pass from  weak to 
limited (and also to  Pettis-norm) convergence  in $\l1$.   
It is shown that tightness is  a 
necessary and sufficient condition to pass from 
limited to strong convergence. 
Other implications between several  modes of convergence 
in $\l1$ are also studied. 
\endabstract 
\endtopmatter 
\document
\baselineskip 14 pt 

\heading 
1. Introduction
\endheading  
\vskip 10pt 

Vaguely speaking,   
a relatively weakly compact set   in 
$\rl1$ is relatively norm compact  if  the functions in 
the set do not oscillate too much.   
Specifically, 
a relatively weakly compact subset of $\rl1$ 
is relatively norm compact if and only if 
it satisfies the Bocce criterion (an oscillation condition) [G1, G2].  
However, the set of constant functions of norm at most one 
in $\l1$ already  shows that (for a reflexive infinite-dimensional 
Banach space $E$), in the Bochner-Lebesgue 
space  $\l1$, more 
care is needed in  order to pass from weak to strong compactness.  
In Section~2,  we   
extend  from $\rl1$ to $\l1$ the above weak-to-norm  result, along with the 
sequential analogue.  
In Section~3, limited convergence (a weakening of strong convergence 
[B1,B2]) is examined.  Limited convergence provides an extension of the 
Lebesgue Dominated Convergence Theorem  to $\l1$.  
Necessary and sufficient conditions to 
pass from  weak to 
limited convergence    
are given.  
In Section~4, the concept of tightness helps to extend the 
results from the previous two sections.  
In Section~5, convergence in the Pettis norm, 
a weakening of strong convergence along lines distinct 
from limited convergence, is examined.  
Similarly, necessary and sufficient conditions to 
pass from  weak to 
Pettis-norm  convergence    
are given.  
In the study,  implications between 
several modes of convergence on $\l1$  
are  examined.  
 
Throughout this paper 
$(E,\|\cdot\|)$  is a Banach space with dual  $E^{*}$ 
and  $B_E$ is the closed unit ball of $E$. 
The triple  $(\Omega, \F ,\mu)$ is a finite measure space.     
Without loss,  we take $\mu$ to be a probability measure.   
For $B\in\F$,  we often examine 
the collection  $\F^+(B)$  of all measurable subsets of $B$  
with (strictly) positive measure and denote $\F^+(\Omega)$ 
by just $\F^+$.      
By $\l1$ we denote the (prequotient) space of all Bochner 
$\mu$-integrable functions from $\Omega$ into $E$.    
On this space the classical $\l1$-seminorm is given by
$\| f \|_{1} := \int_{\Omega} \| f \| d\mu$ and  convergence in this
seminorm is called {\it strong} convergence.

Recall [IT] that the dual of $(\l1,\| \cdot \|_{1})$ is
the (prequotient) space $\dl1$ of 
{\it scalarly} measurable bounded functions from $\Omega$ into $E^{*}$.   
The subspace $\li$
of $\dl1$ consisting of the {\it strongly}
measurable functions actually coincides with $\dl1$
if and only if $E^{*}$ has the {\it Radon Nikodym property} (RNP);
cf. [DU, IT].
Convergence in the corresponding weak topology
$\sigma(\l1, \dl1)$ is called {\it weak} convergence. 
We will also consider the $\sigma(\l1, \li)$-topology on $\l1$.

Also recall that a subset $K$ of $\l1$ 
functions is {\it uniformly integrable} if 
$$
\lim_{c\to\infty} \sup_{f\in K} \int_{[||f|| \ge c]} || f || \, d\mu  =
0  \ . 
$$  
It is  well known  [N] that $K$ is  uniformly integrable if and only if 
it is  bounded 
(i.e. $\sup_{f\in K} ||f||_{\l1}$ is finite)
and 
equi-integrable, i.e. 
$$
\lim_{\mu(A)\to 0^+} \sup_ {f\in K} \int_A || f || \, d\mu  = 0 \ .
$$

All notations and terminology, not otherwise explained, 
are as in [DU, IT, or N].

\heading
2. Weak  vs. Strong Convergence in $\l1$
\endheading 
\vskip 10pt

Our goal is to 
determine precisely when 
(via an oscillation condition)  a weakly convergent sequence   
is also strongly convergent, along with the nonsequential analogue.

For $f \in  \l1$ and $A \in \F $,
the {\it average  value}  
and the {\it Bocce oscillation}   of $f$ 
over $A$ are (respectively) 							
$$
\align
m_A(f) ~&:=~ \frac{\int_A f \, d\mu}{\mu (A)}  \\
\text{Bocce-osc }f \big |_A ~&:=~ 
\frac{\int_A || f- m_A(f)   || \, d\mu }{\mu (A) } 
\endalign
$$
observing the convention that 0/0 is 0. 
The following  elementary inequalities are  useful 
$$
\align
\biggl| 
 \text{Bocce-osc }f \big |_A - \text{Bocce-osc }g \big |_A 
\biggr|  
&\le \text{Bocce-osc }(f-g) \big |_A \\
\text{Bocce-osc }(f+g) \big |_A
&\le
 \text{Bocce-osc }f \big |_A + \text{Bocce-osc }g \big |_A   \\
\mu(A) ~ 
\text{Bocce-osc }f \big |_A  \ 
&\le 2 \int_A || f || \ d\mu  \ . 
\tag2.1
\endalign  
$$
In the spirit of [G1], we consider the following oscillation conditions. 
\definition{Definition 2.1}  [B3]      
A sequence $(f_k )_{k=1}^\infty$ of functions in $\l1$ satisfies the 
{\it sequential Bocce criterion} if   
for each subsequence $( f_{k_j} )$ of $(f_k )$,   
each $\ep$, and  each $B$ in $\F^+$  
there is a set $A$ in $\F^+(B)$   
such that
$\liminf_j \text{Bocce-osc }f_{k_j} \big |_A < \e  \ . $
\enddefinition\flushpar
  
\definition{Definition 2.2}  [G1]   
A subset $K$ of $\l1$ satisfies the  {\it Bocce criterion} if 
for each $\ep$ and  each $B$   
in $\F^+$  
there is a finite collection  $\Cal A$ of sets in $\F^+(B)$    
such that for each $f$ in $K$ 
there is a set  $A$ in $\Cal A$ satisfying  
$ \text{Bocce-osc }f \big |_A < \e  \ . $
\enddefinition\flushpar 

It is known [G1, G2] that a relatively weakly compact subset of $L^1_{\Bbb R}$ 
is relatively norm compact if and only if it satisfies the Bocce criterion. 
We now extend to  $\l1$.  
\proclaim{Theorem 2.3} 
A sequence $ (f_k) $ in  $\l1$   converges 
{\rm strongly}  to $f_0$ in  $\l1$
if and only if  
\roster 
\item $ (f_k) $  converges {\rm weakly} to $f_0$ in  $\l1$
\item $ (f_k )$ satisfies the sequential Bocce criterion 
\item $\Delta_B := \{ m_B (f_k) :k\in N \}$ is relatively norm compact in $E$ 
       for each $B\in\F^+$. 
\endroster  
Condition {\rm  (1)}  may be replaced with 
\roster 
\item"{ (1')}"  $ (f_k) $  converges  to $f_0$ in  $\l1$ in the 
      $\sigma(\l1, \li)$-topology.  
\endroster
Also, condition {\rm  (3)} may be replaced with 
\roster 
\item"{ (3')}"  $\lim_k || m_B(f_k) - m_B (f_0) || = 0 $ 
       for each $B\in\F^+$. 
\endroster
\endproclaim \flushpar  
Note that Theorem~2.3 need not hold if one replaces 
condition (1) (resp\.~(1')\,) with   
$(f_k)$ is Cauchy in the 
weak (resp\. $\sigma(\l1, \li)$-) topology 
since $\l1$ need not be  sequentially  complete in  
this    topology.   
Recall  that $\l1$ is weakly  sequentially  complete 
if and only if $E$ is  [T]; on the other hand,      
$\l1$ is $\sigma(\l1, \li)$-sequentially complete 
if and only if $E$ is weakly sequentially    
complete and has the RNP (cf.  [BH1], [SW]).

There is a set-analogue of Theorem~2.3:   
\proclaim{Theorem 2.4} 
A subset $K$ of $\l1$ is relatively {\rm norm} compact if and only if  
\roster 
\item $K$ is relatively weakly compact 
\item $K$ satisfies the  Bocce criterion 
\item $\Delta_B := \{ m_B (f) : f\in K \}$
      is relatively norm compact  in $E$  for each $B\in\F^+$. 
\endroster  
Condition  {\rm (1)}  may be replaced with 
\roster 
\item"{ (1')}"  $K$ is relatively compact in the 
      $\sigma(\l1, \li)$-topology.  
\endroster
\endproclaim 

Note that  the above  condition (3) is indispensable,    
as shown by Example 3.2         
to come. 
In general, 
if $(f_k )$ is weakly convergent 
(resp. $K$ is relatively weakly compact), 
then the corresponding sets $\Delta_B$ are relatively 
weakly compact in $E$. 
Thus if $E$ is finite-dimensional,  
then condition (3) in the above theorems is not necessary.
 
It is possible to prove Theorems 2.3 and 2.4 by  using 
methods similar to those in [G2].   
Here ideas from 
both [B3] and [G2] are combined.   
The following elementary lemmas are useful.    

\proclaim{Lemma 2.5} 
If  $f $ is in $\l1$, then for each  
$\ep$ and $B\in \F^+$ there is  
a set $A$ in $\F^+(B)$  such that 
$\text{\rm Bocce-osc }f \big |_{A_0} < \e$ for each  subset $A_0$ of $A$. 
\endproclaim 

\demo{Proof} By strong measurability of $f$ in $\l1$ 
and Egorov's Theorem, there exists a 
sequence of simple functions converging almost uniformly 
to $f$.  In combination with (2.1), 
the remainder of the proof is clear.  
\qd 
\enddemo 

\proclaim{Lemma 2.6}  
Let  $\phi \: \Omega \to [0, +\infty]$ be measurable. 
If for each $\ep$ and each  $B$ in $\F^+$ 
there exists  a set $A$ in $\F^+(B)$   
such that $m_A (\phi) < \e$,  
then $\phi(\omega) = 0$ for a.e. $\omega$. 
\endproclaim 

\demo{Proof} 
Fix $\ep$.  Let $B$ be the set of all $\omega \in \Omega$ with 
$\phi (\omega) \ge 2 \e$.  If $B \in \F^+$, then for the 
corresponding  set $A$ in $\F^+(B)$  
we would have $ 2\e \mu(A) < \e \mu(A)$, 
which cannot be.  So $B$ must be a null set.  
\qd 
\enddemo 

\demo{Proof of Theorem 2.3} 
Consider a sequence $ (f_k) $ in  $\l1$   which  converges 
strongly  to $f_0 $. 
Conditions (1)  and (3)  follow immediately. Also, by
(2.1) one has  that
$$  
\mu(A) 
\biggl|
 \text{Bocce-osc } f_{k} \big |_{A}
- \text{Bocce-osc } f_{0} \big |_{A} 
\biggr|
~ \leq ~
2 \int_{\Omega} \| f_{k} - f_{0} \| \, d\mu \rightarrow 0 
$$ 
for each $A$ in ${\Cal F}^{+}$. 
By Lemma 2.5 the singleton 
$\{ f_{0} \}$ satisfies the Bocce criterion.  Thus condition (2) also holds. 

As for sufficiency of (1), (2), and (3),    
note that to prove strong convergence it 
is enough to show that any subsequence $(f_{n})$ of $(f_{k})$ contains 
a further subsequence which converges strongly to $f_{0}$. 
By condition (1)       
the set  $(f_{k})$ is uniformly integrable; hence $(\| f_{k} - f_{0} \|)$ 
must  also be uniformly integrable. So
the subsequence  $(f_{n})$ contains a  further subsequence
$(f_{n_{j}})$ such that $(\| f_{n_{j}} - f_{0} \|)$ converges
weakly to some (nonnegative) function $\phi$ in $\rl1$. 
We shall show that Lemma 2.6 applies to $\phi$; this then 
gives $\phi = 0$ a.e.,  which  finishes the proof.
To show that the lemma applies, first note that by Lemma 2.5
(applied to $f_{0}$) and the given Bocce property (2), the sequence 
$(f_{k} - f_{0})$ also satisfies the  sequential Bocce criterion.
Now fix $\epsilon > 0$ and $B$ in ${\Cal F}^{+}$.
Let $A$ in $\F^+(B)$ be as in Definition 2.1 applied to 
the subsequence $(f_{n_{j}} - f_{0})$ of   
$(f_{k} - f_{0})$, thus     
$$
\liminf_{j} \text{Bocce-osc } (f_{n_{j}} - f_{0}) \big |_{A} < \epsilon.
$$
But by the triangle inequality
$$
\text{Bocce-osc } (f_{n_{j}} - f_{0}) \big |_{A} \geq
\frac{1}{\mu(A)} \int_{A} 
\left[ \| f_{n_{j}} - f_{0} \| - \| m_{A}(f_{n_{j}} - 
f_{0}) \| \right] \, d\mu \ , 
$$
so by weak convergence of $(\| f_{n_{j}} - f_{0}\|)$ to $\phi$ and by 
the given property (3), this leads us to $m_{A}(\phi) < \epsilon$,
which is precisely what  is needed to  apply
Lemma 2.6.  
\qd 
\enddemo

A close look at the proof reveals that 
the conditions may be slightly  weakened.  
Using terminology and results  to come in  Section~3, 
note that 
condition (1) may be replaced with the  two 
conditions that 
$(f_k)$ is  uniformly integrable 
and that 
$(f_k)$ converges scalarly weakly (see Definition~3.3)     
 to $f_0$ in $\l1$.  
These two conditions are equivalent to (1'),   
as noted in Remark~3.7.  
Also, condition (3') is equivalent to 
the two conditions 
that $(f_k)$ converges scalarly weakly to $f_0$ 
and condition (3).  
Thus, under condition (1) or (1'), 
condition (3) is equivalent to (3').     

\demo{Proof of Theorem 2.4}
It is well-known and easy to check that a 
subset  $K$ of $\l1$ is relatively strongly compact 
if and only if it satisfies condition (3) and for each $\eta > 0$ 
there is a 
finite measurable partition $\pi$ of $\Omega$  such that 
$\int_{\Omega} \| f - E_{\pi}(f) \|d\mu < \eta$ for each $f$ in $K$.     
Here $E_{\pi}(f)$ denotes the 
conditional expectation of $f$ relative to the finite algebra generated
by $\pi$. 

Consider a  relatively strongly compact subset $K$ of $\l1$. 
Clearly conditions (1) and (3) are satisfied.   
To see that condition (2) holds, 
fix $\ep$ and $B \in {\Cal F}^{+}$. 
Next, from the above observation, 
find the  partition $\pi := \{ A_{1}, \cdots , A_{N} \}$ corresponding to 
$\eta := \epsilon \mu(B)$.   
Put $\Cal A = \{ A_{i} \cap B \in \F^+ \: A_i \in \pi \}$. 
Fix $f$ in $K$.  
Since 
$$
\split
\sum_{i} \mu(A_{i} \cap B) \text{Bocce-osc } f \big |_{A_{i} \cap B}
&\leq \sum_{i} \mu(A_{i}) \text{Bocce-osc } f \big |_{A_{i}} \\
&= \int_{\Omega} \| f - \sum_{i} m_{A_{i}}(f) 1_{A_{i}} \|d\mu <
\epsilon \mu(B) \ , 
\endsplit
$$
for at least one $A_{i} \cap B \in \Cal A$ we have 
$\text{Bocce-osc } f \big |_{A_{i} \cap B} < \epsilon$.

As for the sufficiency of (1), (2), and (3),   
note that it is enough to show relative
strong  {\it sequential} compactness of $K$.  
So consider a sequence $(f_k)$ in $K$.  
By condition~(1), there is a subsequence 
$(f_{k_j})$ of $(f_k)$ that converges 
weakly to some function $f_0$ in $\l1$ 
while condition~(2) implies that $(f_{k_j})$ satisfies the  
sequential Bocce criterion.  Now an appeal  to Theorem~2.3 shows  that 
$(f_{k_j})$ converges strongly, as needed.

As for replacing  (1) with (1'), 
recall [BH2] that  for the 
$\sigma(\l1, \li)$-topology, 
relatively compact sets and 
relatively  sequentially compact sets 
coincide. 
\qd
\enddemo 

Section~5 gives several variations of the 
Bocce criterion which  also 
provide necessary and sufficient conditions 
to pass from weak to strong convergence (resp. compactness).

\heading  
3. Limited Convergence 
\endheading 
\vskip 10pt 

This section examines limited convergence, 
a weakening of strong convergence [B1].    
Limited convergence  
provides an extension to $\l1$  of the  
Vitali Convergence Theorem (VCT), 
thus also of the  Lebesgue Dominated Convergence 
Theorem (LDCT).   
Furthermore, it extends   the  
previous section's  results.   
In the next section, a tightness condition  
ties    together limited and strong convergence 
and thus extends the results of this section.

Let $\G$ be the collection of all functions 
$g \: \Omega \times E \to \Bbb R$ satisfying 
\roster 
\item"{(i)}" $g(\omega, 0) = 0 $ \  for each  $\omega $ in  $\Omega$ 
\item"{(ii)}" $g(\omega, \cdot)$  \ is weakly $\sigma(E, E^*)$-continuous
for each $\omega$ in $ \Omega$  
\item"{(iii)}" $|g(\omega,\cdot)| \leq C\|\cdot\| + \phi(\omega)$ \  
for each $\omega $ in $\Omega$, for some $C > 0$ and $\phi $ in $ \rl1$
\item"{(iv)}" $g(\cdot,f(\cdot))$ \  is $\F$-measurable for each
$f $ in  $\l1$.
\endroster 
An example of such a function $g$ in $\G$ is given by
$ g(\omega,x) = \sum_{i=1}^n  |x^*_i (x) | \, 1_{A_i} (\omega)$ 
where $x^*_i \in E^*$ and $A_i \in \F$. 
The function $g$ given by $g(\omega,x) = || x ||$ is in 
$\G$ if $E$ is finite-dimensional (for only then does $g$ satisfies (ii)).  
The class $\G$ serves as a ``test class'' for limited convergence   
(see Remark~3.9).        
\definition{Definition 3.1} 
A sequence $(f_{k})$ of functions in $\l1$ 
converges {\it limitedly} to $f_{0}$ in $\l1$ if 
$\lim_{k}\int_\Omega g(\omega, f_k (\omega) - f_0 (\omega)) \ d\mu(\omega) = 0$ 
for each $g\in\G$.  
\enddefinition

Strong convergence implies limited convergence.   
For first note that  a sequence  converges limitedly to $f$ 
if each subsequence has a further subsequence which  
converges limitedly~to~$f$.  
Next note that a strongly convergent sequence    
has the property that 
each subsequence has a further subsequence which is   
pointwise a.e. strongly  convergent.    
Lastly  note that any  uniformly integrable 
sequence $(f_k)$ which is  a.e.  weakly  null  
(i.e.  there is a set $A$ of full measure 
such that if $x^* \in E^*$ and $\omega\in A$ 
then $x^*f_k(\omega)$ converges to zero) converges limitedly.  
To see this, fix $g\in \G$ and put $h_k(\omega) = g(\omega, f_k (\omega))$. 
Condition (iii) gives 
that the set $( h_k )$ is uniformly integrable.  
Conditions (i) and  (ii) give that  $( h_k )$  is a.e.-convergent to 0.  
So $( h_k )$ converges strongly to zero 
and  so $( f_k )$ converges limitedly.   
 
If $E$ is finite-dimensional    
then strong and limited convergence coincide 
(consider  $g \in \G$ given by $g(\omega,x) = || x ||$).  
However, as seen by  modifying the  next example, 
for any infinite-dimensional reflexive space $E$      
there is a sequence of $\l1$ functions which converges  
limitedly but not strongly. 

\example{Example 3.2 (limited $\nRightarrow$ strong)} 
Take $(\Omega, \F,\mu)$ to be the interval $[0,1]$, equipped with
the Lebesgue $\sigma$-algebra and measure and  $E := \ell^{2}$. 
Setting $f_{k}$ identically equal to the $k$-th unit vector $e_{k}$
in $\ell^{2}$ gives a sequence $(f_{k})$ which converges limitedly 
but not strongly  to the null function. 
\endexample  

Limited convergence implies weak convergence since 
for each $b\in\dl1$
the function $g$  defined by 
$g(\omega, x) =  \langle x,b(\omega) \rangle$ is in $\G$. 
As for the converse implication, even  
for finite-dimensional $E$      
weak convergence does not imply limited convergence.  

Towards a variant of the VCT--LDCT for a sequence $( f_k )$ in $\l1$,
we examine  the corresponding sequences 
$(x^*( f_k))$ in $\rl1$  for  $x^*$ in $E^*$.    
\definition{Definition 3.3} 
A sequence $(f_k)$ of functions in $\l1$ converges 
{\it scalarly strongly}  
(resp. {\it scalarly in measure, scalarly weakly\,}) 
to $f_{0}$ in $\l1$ if      
the corresponding  sequence $(x^*(f_{k}))$  in $\rl1$   
converges 
{\it strongly} (resp.  {\it in measure, weakly\,})  to  $x^*(f_{0})$ 
for each $x^*$ in $E^{*}$.  
\enddefinition\flushpar  
Note  the following  chain of {\it strict} implications: 
$$
\text{\it strong} 
\Rightarrow 
\text{\it limited} 
\Rightarrow
\text{\it scalarly strong} 
\Rightarrow
\text{\it scalarly in measure} \ . \tag{3.1}
$$   
Since for $x^*  \in E^*$  functions of the form 
$ g(\omega,\cdot) =  |x^* (\cdot)  | \,  1_\Omega (\omega)$    
are in $\G$,   
limited convergence implies scalarly strong convergence.  
The other implications in (3.1) are clear. 

Furthermore, the implications are strict. 
Example 3.2 showed the first one is not reversible. 
The last implication is not reversible even for $E = \Bbb R$. 
The next example shows that the second implication 
is  also strict. 

\example{Example 3.4 (scalarly strong $\nRightarrow$ limited)} 
Take $(\Omega, \F, \mu)$, $E$, and $(e_k)$  as in Example~3.2. 
Let $I^j_i$ be the dyadic interval $[(i-1)2^{-j},i 2^{-j})$ 
for $j \in {\Bbb N}$ and  
$i = 1, \cdots, 2^{j}$.   
Consider the sequence $(f_{k})$ of the functions
$f_{k} \: [0,1] \rightarrow \ell^{2}$ given by
$f_{k}(\omega) = 1_{I^{k}_{1}}(\omega) 2^{k} e_{k}$.
Since for every $y^* := (y^j)_j$ in $E^* \approx \ell^{2}$
$$
\int_\Omega | y^* \left(f_{k}\left(\omega\right)\right) |\ d\mu(\omega) 
~=~   |y^{k}|  \ ,  
$$ 
$(f_{k})$ converges scalarly strong  to the null function. 
But for the test function 
$g\left(\omega, (x_j) \right) 
= \sum_{j = 1}^{\infty}  x_{j}~1_{I^{j+1}_2}(\omega) $  \ in $\G$ 
$$
\int_\Omega g\left(\omega,f_{k}\left(\omega\right)\right)\ d\mu(\omega) =
\int_{I^k_1} 2^{k} ~ 1_{I^{k+1}_2} \, d\mu = \frac{1}{2} \ .
$$
So $(f_k)$ does not converge limitedly to the null function. 
\endexample \flushpar 
Note that a scalarly strongly convergent sequence need not be
uniformly integrable (as Example~3.4 shows). 
However, a limitedly convergent sequence, being also weakly convergent, 
is  necessarily uniformly integrable.

Limited convergence provides the following extension of the  
VCT--LDCT to $\l1$.  

\proclaim{Theorem 3.5} 
Let $E^*$ have the RNP.   
If a uniformly integrable sequence $(f_k)$ converges 
scalarly in measure to $f_0$ in $\l1$,
then it also converges limitedly to $f_0$.
\endproclaim \flushpar
The necessity of the uniform integrability condition has 
already been noted 
while the    
necessity of $E^*$ having  the RNP 
follows from Remark~5.4.       
The proof of Theorem 3.5 uses the following lemma.

\proclaim{Lemma 3.6} 
A uniformly integrable sequence $(f_k)$ of $\l1$ functions converges limitedly 
to the null function provided that, for each $N \in \Bbb N$, the  
sequence $(f_k^N)_k$ of truncated functions converges limitedly 
to the null function, where 
$f_k^N := f_k \ 1_{[ || f_k || \le N ]}$. 
\endproclaim 

\demo{Proof} 
Fix $g\in\G$ with $|g(\omega,\cdot)| \leq C\|\cdot\| + \phi(\omega)$. 
Now 
$$
\align 
\left|\, \int_\Omega 
\left( g\left(\omega, f_k  \left(\omega\right)\right)  -  
       g\left(\omega, f_k^N\left(\omega\right)\right) 
\right) \ d\mu \, \right| ~&=~  
\left |\, \int_{[ || f_k || > N ]}  g(\omega, f_k(\omega)) \ d\mu \, \right| \\ 
&\le C \int_{[ || f_k || > N ]} || f_k || \ d\mu 
     \ + \ \int_{[ || f_k || > N ]} \phi \ d\mu \ , 
\endalign 
$$ 
so by uniform integrability of $(f_k)$ it follows that
$$
\lim_{N \to \infty} \sup_k 
\left|\, \int_\Omega 
\left( g\left(\omega, f_k  \left(\omega\right)\right)  -  
       g\left(\omega, f_k^N\left(\omega\right)\right) 
\right) \ d\mu \, \right| ~=~ 0  \ . 
$$ 
The lemma now follows with ease. 
\qd
\enddemo

\demo{Proof of Theorem 3.5} 
Without loss of generality, 
we assume that  $f_{0} = 0$ and 
(using the previous lemma) that 
the $f_{k}$'s are  uniformly bounded.  
Note that we may
also assume that $E^{*}$ is separable.
Indeed, by the Pettis measurability theorem [DU, Theorem~II.1.2],   
there is a separable
subspace $E_{0}$ of $E$ such that the $f_{k}$'s are essentially valued
in $E_{0}$. 
Because $E^{*}$  has the RNP, $E_{0}^{*}$
must be separable [DU, Corollary~VII.2.8].
Moreover, if $(f_{k})$ converges limitedly to 0 in
${\Cal L}_{E_0}^{1}$, then it also does so in $\l1$.  

As noted earlier, it is  enough to show
that every subsequence of $(f_{k})$ has a further subsequence that 
is a.e. weakly null.  We  assume  (w.l.o.g.) that this former
subsequence is actually the entire sequence $(f_{k})$.

Now let $(x^{*}_{i})$ be a countable dense subset of $E^{*}$. 
For each $i$
the sequence $(x^{*}_{i} (f_{k}) )$    
converges in measure to zero.
So there exists a subsequence $(f_{k_{j}})$ such that  for a.e. $\omega$ 
$$
\lim_{j}  \  x^{*}_{i} (f_{k_{j}}(\omega))~=~0  \ .  
$$
By a Cantor diagonalization argument  
there is a set  $A$  of full measure  
and a subsequence $(f_{k_{p}})$ such that 
$\lim_{p}  x^{*}_{i} (f_{k_{p}}(\omega)) = 0$
for each fixed $i$ and each $\omega $ in $A$. 
Since the $f_k$'s are uniformly bounded and  
$(x^{*}_{i})$ are dense in $E^*$,   
this pointwise limit property extends so that 
$\lim_{p} \, x^{*}(f_{k_{p}}(\omega)) = 0$
for each fixed $x^{*} $ in  $E^{*}$ and each $\omega $ in $A$. 
Thus, $(f_{k_{p}})$  is a.e. weakly null, as needed.  
\qd 
\enddemo 

Limited convergence also provides an extension of the  
results from the previous section; namely, 
it is possible to pass from weak to limited convergence 
via an oscillation condition. 
The following string of {\it strict}  implications 
summarizes the ideas thus far.  
$$
\text{\it scalarly strong} \Rightarrow \text{\it scalarly weak}
                \Leftarrow  \text{\it $\sigma(\l1, \li)$-topology}
                \Leftarrow \text{\it weak} \ . 
$$

\remark{Remark 3.7} 
A scalarly weakly convergent sequence 
converges in the $\sigma(\l1, \li)$-to\-pol\-o\-gy 
if and only if it is uniformly integrable.  
(Recall that the simple functions 
are not  dense in $ \li$ for infinite-dimensional $E$.) 
Convergence in the $\sigma(\l1, \li)$-topology 
implies weak convergence if and 
only if  $E^*$ has the RNP  [cf. DG] .     
\endremark

In the light of these observations 
and Theorem 3.5, 
we have the following variant of 
Theorem 2.3  for limited convergence. 

\proclaim{Theorem 3.8} 
Let $E^*$ have the RNP. 
A sequence $( f_k) $ of  $\l1$ functions  converges 
{\rm limitedly}  to $f_0$ in $\l1$
if and only if  
\roster 
\item $( f_k) $  converges {\rm weakly} to $f_0$ in $\l1$  
\item $( x^*(f_k) )$ satisfies the 
      sequential Bocce criterion for each $x^*$  in  $E^*$. 
\endroster  
Condition {\rm (1)} may be replaced with 
\roster 
\item"{ (1')}"  $ (f_k) $  converges  to $f_0$  in $\l1$ in the 
      $\sigma(\l1, \li)$-topology.  
\endroster
\endproclaim 

\remark{Remark 3.9}
Limited convergence for separable reflexive $E$ was introduced in 
[B1, B2].  There, the condition (iv) is replaced with  
\roster
\item"{(iv')}" $g$ is $\F \otimes {\Cal B}(E)$-measurable,
\endroster
Of course (iv') always implies (iv). To see that (iv) implies (iv')  
if $E$ is separable, consider a function $g$ which satisfies (iv). 
For each $k\in\Bbb N$, 
write $1_E = \sum_{n} 1_{E^k_n}$ where  
$E^k_n  \in {\Cal B}(E)$ and  the diameter of $E^k_n$ 
is less than    $\tfrac{1}{k}$.  
Choose $x^k_n \in E^k_n$.  Define  
$g_k \: \Omega \times E \to \Bbb R$ by 
$$  
  g_k(\omega,x) =  \sum_n g(\omega, x^k_n) ~1_{E^k_n} (x)  \ .
$$
Since each 
$g_k$ is $\F \times {\Cal B}(E)$-measurable 
and 
$g_k$ converges to $g$ almost everywhere, 
$g$ is also $\F \otimes {\Cal B}(E)$-measurable.  
\endremark

\heading
4.  The Tightness Connection
\endheading 
\vskip 10pt

The concept of tightness  links  strong and limited convergence.  
In this section, we assume that $E$ is a separable Banach space.  
Tightness is considered here with respect to the norm topology 
on $E$ and only for functions.   
The following formulation of tightness is given in [B4].  

\definition{Definition 4.1}  
A subset $L$ of $\l1$ is {\it tight} if there exists an
$ \Cal F\otimes \Cal B(E) $-measurable function 
$h \: \Omega\times E \rightarrow  [0,+\infty]$ such that 
$$
  \sup_{f \in L} 
    \int_{\Omega}  h(\omega,f(\omega)) \,  d\mu(\omega)~<~+\infty  \ 
$$ 
and such that $\{ x \in E : h(\omega,x) \leq \beta \}$ is compact for each 
$\omega \in \Omega$ and each $\beta \in \Bbb R$.
\enddefinition

In [Jaw], the following equivalent formulation of tightness is observed.   

\definition{Definition 4.1$'$}  
A subset $L$ of $\l1$ is {\it tight} if for each $\ep$ 
there exists a measurable multifunction $F_\e$ from $\Omega$ to   
the compact subsets of $E$ such that 
$$
    \mu \left( \{ \omega\in\Omega \: f(\omega) \notin F_\e(\omega) \} \right) 
     \le \e \  
$$ 
for each $f \in L$.   
We say that such a multifunction $F_\e$ is measurable  
(i.e. graph-measurable)  if its graph 
$\{(\omega,x)\in\Omega\times E : x\in F_\e(\omega)\}$
is an $ \Cal F\otimes \Cal B(E) $-measurable subset of $\Omega\times E$. 
\enddefinition

To see the equivalence in one direction, denote the supremum in
Definition 4.1 by $\sigma$  and define $F_\e(\omega)$ as the set of
all $x \in E$ for which $h(\omega,x) \leq \sigma/\e$. 
In the other direction, one obtains a sequence $(F_n)$ of compact-valued
multifunctions by letting $F_n$ correspond to $\e = 3^{-n}$ in  
 Definition~4.1$'$. Without loss of generality
we may suppose that $(F_n(\omega))$ is nondecreasing
(rather than taking finite unions $\cup_{m \leq n} F_m$). 
Now a function $h$   
satisfying the requirements of Definition 4.1 is obtained by setting 
$h(\omega,x) := 2^n$ for $x \in F_{n}(\omega) \backslash
F_{n-1}(\omega)$  with $F_0(\omega) := \emptyset$ 
and $h(\omega,x) :=  +\infty$ for   
$x \in E \setminus \cup_{n} F_n(\omega) $.

In Definition 4.1$'$ we may assume without loss of generality that 
$F_\e(\omega)$ is convex and contains $0$ for each $\omega$ in $\Omega$    
by consider the corresponding multifunction 
$\omega \longmapsto \overline{\text{co}}(F_\e(\omega) \cup \{0\})$.   
The measurability of this new map follows from [CV, Theorem III.40]      
and [HU, Remark (1), p.~163].  
Therefore, if $L$ is tight and $(B_f)_{f\in L}$ is 
a family of sets from  $\F$,  then the set 
$\{ f~1_{B_f} \: f\in L \}$ is also tight. 
Note  that  a bounded sequence in $\l1$ is tight 
if $E$ is finite dimensional (simply take $h(\omega,x) := ||x||$ in 
Definition 4.1). 
For further details on tightness see [B4, B5].  

Recall the following fact [ACV, Th\'eor\`eme~6].
\proclaim{Fact~4.2}
A  uniformly integrable tight subset of $\l1$ is relatively weakly compact. 
\endproclaim \flushpar 
Although weak compactness is not sufficient to 
guarantee that the corresponding subset 
$\Delta_B$ are relatively norm compact 
(consider Example~3.2),  
the following generalization of a result of Castaing [C1]  
shows that    uniform integrability plus tightness 
is sufficient.   

\proclaim{Lemma 4.3}
Let   $L$  be a  tight uniformly integrable  subset of $\l1$.
Then $\Delta_B := \{ m_B (f) :f\in L \}$
is relatively norm compact in $E$
for each $B$ in $\F^+$.
\endproclaim 

\demo{Proof}
Let  the subset $L$  of $\l1$ be  uniformly integrable
and tight.  Since for each $B\in\F^+$
the set $\{ f\ 1_B  \: f\in L \}$
is also uniformly integrable and tight,
it is enough to show that
$\Delta_\Omega$ is relatively norm compact.
Arguing as in Remark (1) on  p.~163 of [HU], we may suppose without
loss of generality that $\F$ is complete.

Fix $\delta >0$.
By the uniform integrability of $L$,
there exist $\alpha >0$ and $\ep$ such that
for each  set
$A$ of measure at most $\e$
we have that
$$
  \sup_{f\in L}
    \int_A   || f || \, d\mu  ~\le~\delta /2
\quad
\text{and}
\quad
\sup_{f\in L}
    \int_{[ || f ||>\alpha]}   || f || \, d\mu
      ~\le~\delta /2
\; .
$$
Let $F_\e$ be a multifunction given by
Definition~4.1$'$  and  $G_\e^\alpha = F_\e\cap\alpha
B_E$ (i.e. $G_\e^\alpha(\omega) =  F_\e(\omega)\cap
\alpha B_E,\,\forall \omega\in\Omega$).
Since $G_\e^\alpha$ is convex compact valued and
integrably bounded (that means
$|| G_\e^\alpha || =
\sup\{ || x || \: x\in G_\e^\alpha(\omega) \}
\in L^1_{\Bbb R_+}$),
the subset
$K_\e^\alpha  = \{\int_\Omega G_\e^\alpha\,d\mu \}$
is convex and compact in $E$ [CV, Theorem~V.15].
Let now $A^\e_f$ be the set of all $\omega \in\Omega$
with $f(\omega) \in F_\e(\omega)$.
Note that $\mu(\Omega\setminus A^\e_f)\le \e$.
Since for each $f\in L$
$$
   \int_{ [ || f || \le \alpha ] }
     f \, 1_{A^\e_f}  \, d\mu \in K_\e^\alpha
\; ,
$$
the set
$\Delta_\Omega^{\e,\alpha} :=
\{ \int_{ [ || f || \le \alpha ] }
f\, 1_{A^\e_f} \, d\mu \: f\in L \}$
is relatively compact in $E$.
Moreover, the distance between
$\Delta_\Omega^{\e,\alpha}$ and
$\Delta_\Omega$ is at most $\delta$ since
$$
|| \  \int_\Omega f \, d\mu -
     \int_{ [ || f || \le \alpha ] } f\, 1_{A^\e_f}
         \,  d\mu \  ||
~\le~
\int_{ [ || f || > \alpha ] }
    || f ||
\, d\mu
 +
\int_\Omega
    || f\, 1_{\Omega\setminus A^\e_f} ||
\, d\mu
~\le~ \delta
$$
for each $f\in L$.  Thus $\Delta_\Omega$ is
relatively compact.
\qd
\enddemo

Measure convergent sequences enjoy tightness.  

\proclaim{Lemma 4.4} 
A sequence in $\l1$ which converges in measure is tight. 
\endproclaim 

\demo{Proof} 
Consider a sequence $(f_k)$ in $\l1$ which converges in measure 
to $f_0$.  
For each natural number $k$, let $\lambda_k $ be    
the bounded non-negative image measure on $E$ induced by $\mu$ and
the measurable function $f_k \: \Omega \rightarrow E$.  
Since $E$ is a Radon  space (thanks to the separability assumption), 
$\lambda_k$ is a Radon (or tight) measure.  For each bounded continuous  
function $\phi \in \Cal C^b(E)$, we have 
$$ 
\lambda_k (\phi) = \int_\Omega \phi(f_k(\omega)) \, d\mu(\omega) \ . \tag{4.1} 
$$ 
It is easy to see that the 
measure convergence of $(f_k)$ to $f_0$ in $\l1$ implies 
the narrow convergence 
(or weak convergence in  the $\sigma(\Cal M^b(E), \Cal C^b(E))$-topology) 
of $(\lambda_k)$ to $\lambda_0$.  
For otherwise there would exist 
$\phi\in\Cal C^b(E)$ and a subsequence 
$(f_{k_j})$ converging almost everywhere to $f_0$ and 
such that $(\lambda_{k_j}(\phi))$ does not 
converge to $\lambda_0(\phi)$. 
But by (4.1), this contradicts the 
Lebesgue Dominated Convergence Theorem. 
Therefore, the sequence $(\lambda_k)$ is tight in $\Cal M^b(E)$ 
in the classical sense 
[S, Appendix Theorem~4],          
which implies that $(f_k)$ is tight   
in the sense of Definitions 4.1 and 4.1$'$.
\qd
\enddemo 
   
From Lemmas~4.3 and ~4.4, 
the following reformulation of 
Theorem~2.3 follows with ease.  

\proclaim{Theorem 4.5}  
A sequence $ (f_k) $ in  $\l1$   converges 
{\rm strongly}  to $f_0$ in  $\l1$
if and only if    
\roster 
\item $ (f_k) $  converges {\rm weakly} to $f_0$ in  $\l1$
\item $ (f_k) $ satisfies the sequential Bocce criterion 
\item $ (f_k) $  is tight. 
\endroster  
\endproclaim \flushpar 

Tightness  connects strong and limited convergence.   

\proclaim{Theorem 4.6}  
A sequence $(f_k)$ of $\l1$ converges strongly
to $f_0$
if and only if
$(f_k)$ is tight and converges limitedly to $f_0$.
\endproclaim 

Before proceeding with the proof of Theorem~4.6, 
we note some immediate corollaries. 

\proclaim{Theorem 3.5 - revisited} 
Let $E^*$ have the RNP and $E$ be separable.   
If a uniformly integrable  tight sequence $(f_k)$ converges 
scalarly in measure to $f_0$ in $\l1$,
then it also converges strongly to $f_0$.
\endproclaim

\proclaim{Theorem 3.8 - revisited} 
Let $E^*$ have the RNP and $E$ be separable. 
A sequence $( f_k) $ of  $\l1$ functions  converges 
{\rm strongly}  to $f_0$ in $\l1$
if and only if  
\roster 
\item $( f_k) $  converges {\rm weakly} to $f_0$ in $\l1$  
\item $( x^*(f_k) )$ satisfies the 
      sequential Bocce criterion for each $x^*$  in  $E^*$   
\item $(f_k)$ is tight.
\endroster  
Condition  \rom{(1)}  may be replaced with 
\roster 
\item"{ (1')}"  $ (f_k) $  converges  to $f_0$  in $\l1$ in the 
      $\sigma(\l1, \li)$-topology.  
\endroster
\endproclaim 
  
The proof of Theorem~4.6 uses the following standard fact 
(compare with Lemma 3.6).     
\proclaim{Fact~4.7}
A uniformly integrable sequence $(f_k)$ of $\l1$ functions converges strongly
to the null function provided that, for each $N \in \Bbb N$, the
sequence $(f_k^N)_k$ of truncated functions converges strongly
to the null function, where
$f_k^N := f_k \ 1_{[ || f_k || \le N ]}$.
\endproclaim 
 
\demo{Proof of Theorem~4.6} 
The implication in one direction follows from our previous work. 
As for the other direction,   
let $(f_k)$ be a tight sequence in $\l1$ which
converges limitedly to $f_0$.    
Because the image measure of $\mu$ under 
$f_0$ is a Radon measure on $E$, the singleton $\{ f_0 \}$ must be
tight. Since the union of two tight sets is again tight, we have that
the set $\{ f_k : k \in {\Bbb N} \cup \{ 0 \} \}$ is also tight; let $h$
correspond to this set as in Definition~4.1.
Without loss of generality, we assume that $f_0$ is the null
function and that the $f_k$'s are uniformly bounded (in 
${\Cal L}^\infty_E$) by some $M > 0$.     
To avoid the non-metrizability of the $\sigma( E , E^\ast )$-topology,  
we use ideas from [B5]. 
By well-known facts about Suslin spaces [S, Corollary~2 of Theorem~II.10],      
there exists a metric 
$d$ on $E$ defining a topology $ \tau_d $ weaker than the weak topology
$ \sigma (E, E^\ast) $ and such that $ (E, \tau_d) $ is a Suslin space.
Define $\phi \: \Omega \times E \rightarrow \Bbb R$ by 
$\phi(\omega,x) := \max(-\| x  \|, -M)$.
For each $\e > 0$,   consider the function 
$\phi_{\e} : \Omega \times E \rightarrow \Bbb R $ given by
$$
\phi_{\e}(\omega,x) := \phi(\omega,x) + \e h(\omega,x).
$$   
From the inf-compactness property of $h$ (see Definition 4.1) it
follows that $\phi_{\e}(\omega,\cdot)$ is also inf-compact on $E$ 
for each $\omega \in \Omega$ and $\e > 0$; in turn, this implies
inf-compactness of the same functions for the weak topology
$ \sigma (E, E^\ast) $ and hence for $\tau_d$.
Moreover, the $\Cal F\otimes\Cal B(E)$-measurability\footnote{For 
any of the three topologies $E$ is a Suslin space; hence, it has
the same Borel $\sigma$-algebra $ \Cal B(E) $.}
of $\phi_{\e}$ is evident.

For each $ \e >0$ and $ p\in\NN $ we define the       
approximate         
function $ \phi_{\e}^p \: \Omega\times E \rightarrow \Bbb R $ \, by
$$
\phi_{\e}^p(\omega,x)=\inf_{y\in E} \{ \phi_{\e}(\omega,y) + p\,d(x,y)\}.
$$
Evidently, for each $\e >0$ the sequence $(\phi_\e^p)$ is (pointwise)
nondecreasing. It is well-known [B4,V1]  that $ \phi_{\e}^p $ has an
$ \Cal F\otimes \Cal B(E) $-measurable modification $\psi_{\e}^p$
(i.e., $\psi_{\e}^p(\omega,\cdot) = \phi_{\e}^p(\omega,\cdot)$ a.e.)
such that for each
$\omega\in\Omega$ the function $\psi_{\e}^p(\omega,\cdot)$
is $ d $-Lipschitz continuous on  $E$ and therefore  is 
$\sigma( E,E^\ast)$-continuous.
Furthermore, as a well-known property of this approximation,
by $\tau_d$-lower semicontinuity 
and boundedness below of $\phi_\e(\omega,\cdot)$,   
we have   
$$
\phi_\e(\omega,x) = \lim_p \uparrow \psi_\e^p(\omega,x) 
$$
for a.e. $\omega$  and each $x \in E$. We  now set    
$\widehat\psi_\e^p(\omega,x) = 
\min(\psi_\e^p(\omega,x) - \psi_\e^p(\omega,0),p)$.
Note that
$ -M - \epsilon h(\omega,0) \leq \widehat{\psi}^p_{\e}(\omega,\cdot) \leq p$ 
for a.e. $\omega$, where $\omega \mapsto h(\omega,0)$ is integrable
in view of Definition 4.1 and $f_0 \equiv 0$. 
For each $\e >0$ and $p \in \NN$, the function $ \widehat\psi_\e^p $
satisfies the conditions (i) to (iv) for the test functions
of $\Cal G$.
Therefore, by the limited convergence of the $f_k$'s, we have
$$
\lim_{k\to\infty}\int_\Omega \widehat\psi_\e^p(\omega,
f_k(\omega)) \, d\mu(\omega)
=0.
$$
It follows that
$$ 
\liminf_{k\to\infty}\int_\Omega \psi_\e^p(\omega, f_k(\omega)) \, d\mu(\omega)
\geq 
\int_\Omega \psi_\e^p(\omega,0) \, d\mu(\omega) \  .
$$
Thus, for each $\e > 0$ and  $p \in\NN$
$$
\liminf_{k\to\infty}\int_\Omega \phi_\e(\omega, f_k(\omega)) \, d\mu(\omega)
  \geq
     \liminf_{k\to\infty}\int_\Omega 
     \psi_\e^p(\omega, f_k(\omega))\, d\mu(\omega)
  \geq
     \int_\Omega \psi_\e^p(\omega,0)\, d\mu(\omega) \ . 
$$
The monotone convergence theorem gives,  
for each $\e > 0$
$$
\align
\alpha_\e 
:=
\liminf_{k\to\infty}\int_\Omega \phi_\e(\omega,f_k(\omega))  \, d\mu(\omega)
&~\geq~
\lim_p \uparrow \int_\Omega \psi_\e^p(\omega,0)\, d\mu(\omega)  \\
&~=~
\int_\Omega \phi_\e(\omega,0)\, d\mu(\omega)
~=~ 
\e \int_\Omega h(\omega,0)d\mu(\omega) \ ,
\endalign
$$ 
thus 
$$
0 \leq \alpha_\e \leq \liminf_{k\to\infty} \int_\Omega \phi(\omega,f_k(\omega))
 \, d\mu(\omega)
+ \e~\sup_k \int_\Omega h(\omega,f_k(\omega))\, d\mu(\omega) \ . 
$$
Since $\phi(\omega,f_{k}(\omega)) = - || f_{k}(\omega)  ||$,
by our
initial assumption, the proof is finished by letting $\e$ go to zero. 
\qd
\enddemo

Fact~4.2 and Theorem~4.5 
gives that a  uniformly integrable tight sequence in $\l1$ which
satisfies the sequential Bocce criterion has  a strongly convergent 
subsequence.   
Recall that a sequence $ (f_k) $ is said to be
bounded if $\sup_k || f_k ||_{\l1} $ is finite.
In the above,  
if we relax  
uniform integrability  to  boundedness, 
we need not have strong subsequential convergence 
(just consider the  sequence 
$(  n\,1_{[0, 1/n]} )_n$ in $\rl1$)
but we do have measure subsequential convergence.   
We can state this result as a strong Biting lemma.

\proclaim{Theorem 4.8} 
Let   $ (f_k) $ be a  bounded   tight sequence  in  $\l1$ 
satisfying the sequential  Bocce criterion.    
Then there exist a subsequence, say $(f_n)$, of  $(f_k)$  
and an increasing sequence $(A_n)$ in $\F$ 
such that 
\roster 
\item $\lim_{n\to\infty} \, \mu(A_n) =  \mu(\Omega)$ 
\item the sequence $(f_n\, 1_{A_n} )$  converges strongly in $\l1$ 
\item the sequence $(f_n \, 1_{\Omega\setminus A_n} )$   
      converges to $0$ in measure.
\endroster  
Therefore,    
the subsequence $(f_n)$ converges in measure.  
\endproclaim
The proof uses Gaposhkin's Biting lemma [Ga, Lemma C], 
which is also referred to as 
Slaby's Biting lemma [cf. C2].

\proclaim{Biting lemma}  
Let   $ (f_k) $ be a   bounded  sequence  in  $\l1$.  
Then there exist a subsequence, say $(f_n)$, of  $(f_k)$  
and an increasing sequence $(A_n)$ in $\F$ 
such that 
\roster 
\item $\lim_{n\to\infty} \, \mu(A_n) =  \mu(\Omega)$ 
\item the sequence $(f_n \, 1_{A_n} )$ is uniformly integrable in $\l1$ .  
\endroster 
Note that (1) implies that
 the sequence $(f_n \, 1_{\Omega\setminus A_n} )$ converges to $0$  in 
       measure.
\endproclaim
 
\demo{Proof of Theorem 4.8} 
Consider a   bounded tight sequence $ (f_k) $ in  $\l1$ which 
satisfies the sequential Bocce criterion.    
Apply the  Biting lemma to find the corresponding 
subsequence, say $(f_n)$, of  $(f_k)$ and sequence $(A_n)$ in $\F$. 
Since $(f_n \, 1_{A_n} )$ is  
uniformly integrable and tight, 
it is relatively  weakly sequentially compact. 
By passing to a further subsequence we can assume that $(f_n \, 1_{A_n} )$ 
converges weakly in $\l1$.  
Since $(f_k)$ satisfies the sequential Bocce criterion, 
using condition (1) it is easy to check that 
$(f_n \, 1_{A_n} )$  also satisfies the sequential Bocce criterion 
(in the definition, for a fixed $B\in\F^+(\Omega)$, apply the criterion 
to $B_0 := B\cap A_N$ for  a sufficiently large $N$).  
Theorem~4.5       
gives that $(f_n \, 1_{A_n} )$ 
converges strongly. 
\qd
\enddemo

\heading 
5. Pettis Norm
\endheading
\vskip 10pt 

This section examines Pettis norm convergence in light of 
the previous sections. 

\definition{Definition 5.1} 
A strongly measurable function $f \: \Omega \rightarrow E$ is 
{\it Pettis integrable} if $x^{*}(f)$ belongs to $\rl1$
for every $x^{*}$ in $E^{*}$ and if for every $B$ in $\F$ there
exists $x_{B}$ in $ E$ such that
$$
\int_{B} x^*(f)d\mu = x^*(x_{B}) \text{ \quad for all }
x^{*} \in E^{*}.
$$
The space   $\p1$ of   
(equivalence classes of) all strongly measurable Pettis integrable functions   
forms a  normed linear space under the 
{\it Pettis} ({\it semi\,}){\it norm}  
$$
  || f ||_{\pn} = \sup_{x^*\in B_{E^*}} \int_\Omega | x^*(f) |  \, d\mu \  . 
$$
Clearly $\p1$ contains $\l1$, to
which we restrict considerations. 
\enddefinition

In general,  
Pettis norm convergence on $\l1$ is incomparable 
with limited convergence but is comparable with the other modes 
of convergence in  chain (3.1).  A parallel chain of 
{\it strict} implications is 
$$
\text {\it strong} 
\Rightarrow 
\text {\it Pettis} 
\Rightarrow
\text {\it scalarly strong}  \ . \tag{5.1}
$$ 
Note that when $E$ is finite-dimensional, the two chains (3.1) and (5.1)   
merge into
$$
\text {\it strong} 
\Leftrightarrow 
\text {\it Pettis} 
\Leftrightarrow 
\text {\it limited} 
\Leftrightarrow 
\text {\it scalarly strong}  \ . 
$$
The implications in  chain (5.1) are clear; 
the following two examples show that  they  
are strict.

\example{Example 5.2 (scalarly strong $\nRightarrow$ Pettis)}    
Example~3.2 suffices here but, for later use,  we consider the following 
variation.     
Take $(\Omega, \F, \mu)$, $E:= \ell^{2}$, and $(e_k)$  as in Example~3.2.  
Consider the  Rademacher-type functions $f_k \: [0,1] \to \ell^2$ 
defined by $f_k (\omega) := e_k r_k(\omega)$ 
where $r_k$ is the $k$-th Rademacher function. 
Clearly, $(f_k)$ converges scalarly strong to the null function  
yet the Pettis norm of each $f_k$ is one.  
\endexample

\example{Example 5.3 (Pettis $\nRightarrow$ strong)}    [P]  
Take $(\Omega, \F, \mu)$, $E$, $(e_k)$, and $(I^j_i)$  
as in Example~3.4.    Consider the sequence    
$(f_k )$ of the integrable functions $f_k \: \Omega \to \ell^2$   
given by
$f_k  (\omega) := \sum_{i=1}^{2^k}  1_{I^k_i} (\omega) e_{ 2^k + i}$.
To see that $(f_k)$ converges in the Pettis norm to the null function, 
fix $y^* := ( y^i)_i  \in B_{\ell^2}$. Put $\overline{y}^*  := ( |y^i| )_i$
and note that 
$$
\align
\int_\Omega | y^*( f_k ) | \, d\mu   
~&=~ \sum_{i=1}^{2^k} |y^{2^k +i}| \  \mu(I^k_i) 
~=~ 2^{-k}  ~\overline{y}^* \left( \sum_{i=1}^{2^k}  e_{ 2^k + i}  \right) \\
~&\le~ 2^{-k}~  \bigl\| \sum_{i=1}^{2^k}  e_{ 2^k + i} \bigr\|_{\ell^2} 
~=~ 2^{-\frac{k}{2}} \ .
\endalign
$$ 
Thus $|| f_k ||_{\pn} \to 0$. 
But  $(f_k)$ does not converge strongly 
since 
$\int_\Omega || f_k ||_{\ell^2} \, d\mu    
~=~ 1$.
\endexample

Example~5.3 illustrates (consider $g_k := 2^{\frac{k}{4}} f_k$) 
that a Pettis-norm convergent sequence need not be  
uniformly integrable.    
Example 3.2 shows that  a  limitedly convergent 
sequence need not converge in the Pettis norm.   
Theorem~3.5 gives that if $E^*$ has the RNP, then  a 
uniformly integrable Pettis-norm convergent sequence  in $\l1$ 
also converges limitedly. 
The following remark shows the necessity of $E^*$ having the RNP.  
\remark{Remark 5.4 \rom{[DG]}} 
A uniformly integrable Pettis-norm convergent sequence also 
converges in the $\sigma(\l1, \li)$-topology and, 
if furthermore $E^*$ has the RNP,  then also weakly. 
But if $E^*$ fails the RNP, then  there is 
an essentially bounded    
sequence which converges in 
the Pettis norm but not weakly (thus not limitedly).
\endremark        
In the case that  $E= \ell^1$, 
this sequence is easy to construct.
\example{Example 5.5 (Pettis $\nRightarrow$ limited)}    
Let $(\Omega, \F, \mu )$ be as in Example 3.2 and let $E = \ell^1$. 
Consider the sequence $( f_k)$ in $\l1$ given by 
$f_k (\omega)  := \frac{1}{k} \sum_{i=1}^k r_i(\omega) e_i$, 
where $e_i$ is the $i$-th unit vector in $\ell^1$ and $r_i$ is the 
$i$-th Rademacher function.  
Note that $( f_k)$ is essentially bounded. 
As for the Pettis norm of $f_k$,  
fix $y^* = (y^i)_i \in E^* = \ell^\infty$.  Since 
$$
\int_\Omega |  y^*(f_k) |  \, d\mu ~=~ 
\frac{1}{k} \int_\Omega
 \left| \sum_{i=1}^k  y^i \, r_i(\omega) \right|  \, \mu(d\omega) \ , 
$$ 
Khintchine's inequality  [cf. D1]     
shows that $||f_k ||_{\pn}$ behaves like 
$\frac{1}{{\sqrt{k}}}$ and so   $||f_k ||_{\pn} \to 0$.   
Thus $(f_k)$ converges scalarly weakly to the null function   
and so if it also  converges limitedly or weakly, 
it does so to the null function.  
But consider  $b\in\dl1\approx{\l1}^*$ given by 
$b(\omega) := (1_{[ r_i = 1]}(\omega))_i$\,, along with the corresponding  
test function  
$g(\omega,x) :=  \langle x, b(\omega) \rangle$.
Since 
$\int_\Omega \langle f_k (\omega), b(\omega) \rangle \ d\mu(\omega)~=~\frac12$ 
we see that  
$(f_k)$ does not  converge limitedly nor  weakly. 
\endexample

At this time there is no analogue to Theorem 3.5 which would allow  
one to pass from scalarly in measure  convergence 
to Pettis-norm convergence when $E^*$ has the RNP.      
Note that if the sequence $(f_k )$ is 
Cauchy in the Pettis norm, then the corresponding subsets 
$\Delta_B$ of $E$ are  relatively norm compact
for each $B\in\F^+$.  
But even for an essentially 
bounded (thus uniformly integrable) sequence $( f_k )$ 
for which the $\Delta_B$ are 
all relatively norm compact, 
the implication {\it scalarly in measure $\Rightarrow$ Pettis} 
does not hold in general, 
as shown by Example~5.2.  

It is possible in certain situations to pass from weak to Pettis-norm
convergence. For this,  a measurement
of the oscillation relative to the Pettis norm is needed.

\definition{Definition 5.6}
For $f \in \l1$ and $A \in \F$ the {\it Pettis Bocce oscillation}
of $f$ over $A$ is
$$
\text{Pettis-Bocce-osc }f \big |_A ~:=~ 
\sup_{x^* \in B_{E^*}} \text{Bocce-osc } x^*(f) \big |_{A}  \ .  
$$
Since       
Bocce-osc $ x^* f \big |_{A} $ is at most
$ \| x^* \| $ Bocce-osc $  f \big |_{A}$, 
the 
Pettis-Bocce-osc $f \big |_{A} $
is at most   
Bocce-osc $f \big |_{A}$.    
\enddefinition
\definition{Definition 5.7} 
A sequence $(f_k )$ of functions in $\l1$ satisfies the 
{\it sequential Pettis Bocce criterion} if   
for each subsequence $( f_{k_j} )$ of $(f_k )$,   
each $\ep$, and  each $B$ in $\F^+$,    
there is a set $A$ in $\F^+(B)$  such that
$\liminf_j \text{Pettis-Bocce-osc }f_{k_j} \big |_A < \e  \ . $
\enddefinition
\definition{Definition 5.8} 
A subset $K$ of $\l1$ is {\it Pettis uniformly integrable}
if the corresponding subset 
$\widetilde K := \{  x^\ast f \: x^\ast\in B_{E^*}, f\in K \}$ 
of $\rl1$ is uniformly integrable. 
\enddefinition   \flushpar 
Clearly, $K$ is Pettis uniformly integrable 
if and only if   
it  is Pettis-norm bounded and 
the corresponding set $\widetilde K$ 
is equi-integrable.

The following variants of Lemma~2.5 and Lemma~2.6, respectfully, 
are useful.   
\proclaim{Lemma 5.9}
The sequential Pettis Bocce criterion is translation invariant.
\endproclaim

\demo{Proof} Let the sequence $(f_k)$ satisfy the sequential Pettis  
Bocce criterion  and fix $f \in \l1$. The fact that $(f_k + f)$
also satisfies the Pettis Bocce criterion follows directly from the
definition, Lemma 2.5, and the observation that (cf. inequalities 2.1) 
$$
\text{Pettis-Bocce-osc } (f_{k} + f) \big |_{A} \le
\text{Pettis-Bocce-osc } f_{k} \big |_{A} +
\text{Pettis-Bocce-osc } f \big |_{A}.
\TagsOnRight\tag"\qed"\TagsOnLeft$$
\enddemo

\proclaim{Lemma 5.10}
Let $(f_k)$ be a Pettis uniformly integrable sequence in $\l1$. If for each
subsequence $(f_{k_j})$ of $(f_k)$, each $\ep$, and  each $B$ in $\F^+$, 
there exists a subset $A$ in $\F^+(B)$ such that
$$
\liminf_j \sup_{x^* \in B_{E^*}} \frac{\int_{A} |x^*(f_{k_j})|d\mu}{\mu(A)}
< \epsilon
$$
then $(f_k)$ converges to $0$ in the Pettis norm.
\endproclaim

\demo{Proof}
Assume $(f_k)$ is  Pettis uniformly integrable  but does not converge to 0 in
the Pettis norm.    
Since $(f_k)$ is  Pettis uniformly integrable,   
the subset $\{ |x^*(f_k)| : x^* \in B_{E^*}, k \in {\Bbb N} \}$ of $\rl1$ is  
relatively weakly compact. So there exists $\ep$, a subsequence
$(f_{k_j})$ of $(f_k)$, a sequence $(x^*_{k_j})$ in $B_{E^*}$,
and $g $ in $\rl1$ such that 
$2 \epsilon < \int_{\Omega} |x^*_{k_j}(f_{k_j})|d\mu$ and
$|x^*_{k_j}(f_{k_j})| \rightarrow g$ weakly in $\rl1$.
Since $2 \epsilon \le \int_{\Omega} g \ d\mu$, the set $B := [g >
\epsilon]$ is in $\F^+$.
For any subset $A$ of $B$ with positive measure
$$
\liminf_j \sup_{x^* \in B_{E^*}} \frac{\int_{A} |x^*(f_{k_j})|d\mu}{\mu(A)}
\ge \liminf_j \frac{\int_{A} |x^*_{k_j}(f_{k_j})|d\mu}{\mu(A)} = 
\frac{\int_{A} g\ d\mu}{\mu(A)} > \epsilon.
$$
Thus the lemma holds.
\qd
\enddemo

The Pettis-norm analogue to Theorems~2.3 and ~3.8 now follows with ease. 

\proclaim{Theorem 5.11} 
A sequence $(f_k)$ in $\l1$ converges in the Pettis
norm to $f_0$ in $\l1$ if and only if
\roster
\item $(f_k)$ is Pettis uniformly integrable 
\item $(f_k)$ satisfies the sequential Pettis Bocce criterion
\item $\lim_k\| m_B(f_k) - m_B(f_0) \| =0$ for each
$B\in\F^+$.
\endroster
\endproclaim

\demo{Proof}
Consider a sequence $(f_k)$ that converges  
in the Pettis norm to $f_0$ in $\l1$.  
It is easy to check that conditions (1) and (3) hold.   
Since for $x^* $ in $ B_{E^*}$ and $A$ in $\F^+$
$$
\biggl| \text{Bocce-osc } x^*(f_k) \big |_{A} 
    - \text{Bocce-osc } x^*(f_{0})  \big |_{A} \biggr| 
~\le~  \frac{2}{\mu(A)} \| f_k - f_0 \|_{\pn}
$$ 
and $\text{Bocce-osc } x^*(f_{0}) \big |_{A} \le \text{Bocce-osc } f_{0}
\big |_{A}$, 
from Lemma 2.5 we see that $(f_k)$ satisfies the sequential Pettis Bocce
criterion.

As for the other implication, consider a    
sequence $(f_k)$ which satisfies
conditions (1), (2) and (3). 
To show that $f_{k} \rightarrow f_{0}$ in the Pettis norm,
we will show that $(f_k - f_0)$ satisfies the conditions of Lemma 5.10.
First note that condition (1) gives that $(f_k - f_0)$ is 
Pettis  uniformly integrable.
Fix  $\ep$  and $B$ in $\F^+$. Consider
a subsequence $(f_{k_j})$ of $(f_k)$. Since $(f_k - f_0)$ satisfies the
sequential Pettis Bocce criterion, there is a set $A$  in $\F^+(B)$     
such that $\liminf_j$ Pettis-Bocce-osc $(f_{k_j} -
f_0) \big |_{A} < \epsilon.$ Since
$$
\sup_{x^* \in B_{E^*}} \frac{ \int_{A} |x^*(f_{k_j} -
f_0)|d\mu}{\mu(A)}
~-~ \| m_{A}(f_{k_j} - f_0) \| \le \text{Pettis-Bocce-osc } (f_{k_j} -
f_0) \big |_{A},
$$
using (3) we see that
$$
\liminf_j \sup_{x^* \in B_{E^*}}  \frac{ \int_{A} |x^*(f_{k_j} -
f_0)|d\mu}{\mu(A)} < \epsilon
$$
as needed. Thus $f_k \rightarrow f_0$ in the Pettis norm.
\qd
\enddemo

Remark~5.4 ties weak convergence into Theorem~5.11.  
\proclaim{Corollary 5.12} 
A sequence $(f_k)$ in $\l1$ converges in the Pettis
norm to $f_0$ in $\l1$  and is uniformly integrable 
if and only if
\roster 
\item $ (f_k) $  converges to $f_0$ in the $\sigma(\l1,\li)$-topology
\item $ (f_k )$ satisfies the sequential Pettis Bocce criterion 
\item $\Delta_B := \{ m_B (f_k) :k\in N \}$ is relatively norm compact in $E$ 
       for each $B$ in $\F^+$. 
\endroster  
Condition {\rm (3)} may be replaced by 
\roster 
\item"{ (3')}"  $\lim_k || m_B(f_k) - m_B (f_0) || = 0 $ 
       for each $B$ in $\F^+$. 
\endroster
Furthermore, if $E^*$ has the RNP, then {\rm (1)} is equivalent to
\roster
\item"{ (1')}" $ (f_k) $  converges to $f_0$ weakly in $\l1$.
\endroster
\endproclaim \flushpar
Note that under (1), conditions (3) and (3') are equivalent.   

Since a Pettis convergent sequence need not be tight (consider  
Example~5.5 along with Fact~4.2), there is no  Pettis-analogue to 
Theorem~4.5.

\heading 
6. Variation of the Bocce Criterion
\endheading 
\vskip 10pt

As noted in this section, 
several variations of the sequential Bocce criterion 
also provided necessary and sufficient conditions to 
pass from weak to strong convergence.   
For a sequence $(f_k )$ of functions in $\l1$, 
consider the following  Bocce-like oscillation conditions.

The sequence $(f_k )$ satisfies 
oscillation condition (B0)  if 
for each $\ep$  
and each $B$ in $\F^+$ 
there is a  set $C$ in $\F^+(B)$  and $N\in\Bbb N$ 
such that 
$$ 
        \text{Bocce-osc }f_{k} \big |_C < \e 
$$ 
for  each  $k\ge N$.    

The sequence $(f_k )$ satisfies 
oscillation condition (B1)  if 
for each $\ep$ there is a finite measurable 
partition $\pi = (A_i)_{i=0}^p$ of $\Omega$ with 
$\mu(A_0)<\e$ 
and $N\in\Bbb N$   
such that 
$$
    \text{Bocce-osc }f_{k} \big |_{A_i} < \e
$$ 
for  each  $k\ge N$ and  $1\le i \le p$.   

The sequence $(f_k )$ satisfies 
oscillation condition (B2)  if 
for each $\ep$ there is a finite measurable 
partition $\pi = (A_i)_{i=0}^p$ of $\Omega$ with 
$\mu(A_0)<\e$ such that    
for each collection $(B_i)_{i=1}^p$ of sets with 
$B_i $ in $\F^+(A_i)$ 
there is $N\in\Bbb N$ such that 
$$
    \text{Bocce-osc }f_{k} \big |_{B_i} < \e
$$ 
for  each  $k\ge N$ and  $1\le i \le p$.

The above 3 oscillation conditions   have appeared in the literature 
[V2, B3, J]  under various names.  
In [J], it is shown that 
$$
   \text{(B2)} \Rightarrow \text{(B1)} \Leftrightarrow \text{(B0)} \ . 
$$ 
The proof that (B2) implies (B0)   
and the proof that (B1) implies (B0) are both  straightforward  
while  the proof that (B0) implies (B1)  involves 
an exhaustion argument. 
It is straightforward [cf. G1]   
to show that (B1) implies the sequential Bocce criterion. 

If the 
sequence $ (f_k) $ in  $\l1$   converges 
strongly  then it satisfies (B2).   
This follows from  minor variations of 
earlier arguments and  noting that Lemma~2.5 may be 
strengthened. 
\proclaim{Lemma 2.5 -- revisited} 
Let $f$ be in $\l1$. 
For every $\ep$  
there is  a finite measurable 
partition $\pi = (A_i)_{i=0}^p$ of $\Omega$ with 
$\mu(A_0)<\e$ such that  
for each collection $(B_i)_{i=1}^p$ of sets with 
$B_i$ in $\F^+(A_i)$   
$$
    \text{\rm Bocce-osc }f \big |_{B_i} < \e   
$$ 
for each $1\le i \le p$.  
\endproclaim  \flushpar  
Thus in Theorem~2.3  (and thus also 
in the related theorems) oscillation    
condition (2) may be replaced with 
the condition that $ (f_k) $ 
satisfies either  oscillation  condition (B2),  (B1),  or (B0).

As for the subset analogue, recall [G1] that 
a subset $K$ of $\l1 $ is a 
{\rm set of small Bocce oscillation} if 
for each $\ep$  
there is a finite measurable partition $\pi= (A_i)_{i=1}^p$ of $\Omega$
such that for each $f$ in $K$ 
$$
      \sum_{i=1}^p \mu(A_i) \text{  Bocce-osc }f\big|_{A_i} \ < \ \e  \ .    
$$
As in the $\rl1$ case [G1], 
a relatively strongly compact set is a set of small 
Bocce oscillation and 
a set of small Bocce oscillation satisfies the Bocce criterion. 
Thus in  
Theorem~2.4  the  oscillation   condition (2) may be replaced by   
the condition that $K$ 
be a set of small Bocce oscillation.

\heading 
Acknowledgment
\endheading 
\vskip 10pt

The authors  express our appreciation 
to the Department of Mathematics at the  University
of South Carolina for providing the financial support which enabled us
to start this work.   
We are grateful  
to M. Valadier  and M. Moussaoui 
for their suggestions and discussions 
concerning the Section~4   
and 
for the anonymous referee's  
careful suggestions that helped clarify many fine details.  
B.~Dawson [Da] has independently obtained results in 
the spirit of this paper's Theorem~2.4.


\enddocument